\documentclass{article}
\usepackage{amsrefs,amssymb,amsmath,amsthm,amsfonts,epsfig,graphicx}
\usepackage[utf8]{inputenc}
\usepackage[T1]{fontenc}
\usepackage[english]{babel}
\usepackage{hyperref}
\hypersetup{
    colorlinks,
    citecolor=black,
    filecolor=black,
    linkcolor=black,
    urlcolor=black
}

\theoremstyle{plain}
\newtheorem{teo}{Theorem}[section]
\newtheorem{thm}[teo]{Theorem}

\newtheorem{lem}[teo]{Lemma}
\newtheorem{prop}[teo]{Proposition}

\theoremstyle{definition}
\newtheorem{df}[teo]{Definition}
\newtheorem{exa}[teo]{Example}
\newtheorem{rmk}[teo]{Remark}

\DeclareMathOperator{\len}{len}

\DeclareMathOperator{\diam}{diam}
\DeclareMathOperator{\clos}{clos}

\DeclareMathOperator{\dist}{dist}

\newcommand{\R}    {\mathbb R}

\newcommand{\Z}  {\mathbb Z}

\renewcommand{\epsilon}{\varepsilon}

\newcommand{\continua}{{\mathcal C}}

\author{Alfonso Artigue}

\title{Pseudo-Anosov maps and continuum theory}
\begin{document}
\date{\today}
\maketitle

\begin{abstract}
In the hyperspace of subcontinua of a compact surface 
we consider a second order Hausdorff distance. 
This metric space is compactified in such a way that 
the stable foliation of a pseudo-Anosov map is naturally identified with a hypercontinuum. 
We show that negative iterates of a stable arc converges to this hypercontinuum in the considered metric.
% 
% We apply tools from continuum theory for modelling the stable and unstable foliations 
% of pseudo-Anosov maps. We formalize the idea that negative iterates of a stable arc converges to the 
% stable foliation.
Some dynamical properties of pseudo-Anosov maps, as topological mixing and the density of stable leaves, are generalized for cw-expansive homeomorphisms 
of pseudo-Anosov type
on compact metric spaces.
\end{abstract}
% \tableofcontents
\section{Introduction}
In order to motivate some of the ideas in this paper let us start with an example.
Let $f\colon T^2\to T^2$ be the Anosov diffeomorphism of the two dimensional torus defined by 
$f(x,y)=(2x+y,x+y)$. 
The local stable sets are arcs and they form the stable foliation of $f$.
Take $A^s\subset T^2$ a stable arc. In this case 
% $\diam(f^n(A^s))\to 0$ as $n\to+\infty$ and 
$f^n(A^s)$ converges to $T^2$ in the Hausdorff metric as $n\to-\infty$. 
% This property was investigated in \cite{Ka2} in relation with the concept of \emph{full expansivity}.
The Hausdorff metric $\dist_H$ between compact subsets is recalled in Equation (\ref{eqHdist}) below.
Reversing time we obtain a similar result for an unstable arc $A^u$. 
Then, $f^{-n}(A^s)$ and $f^n(A^u)$ have a common limit, the torus, as $n\to+\infty$ in the Hausdorff metric. 
This is what we find if we consider limits of stable and unstable arcs in the Hausdorff metric,
but these limits do not show that the arcs are \emph{winding} around the surface in \emph{different directions}.
The purpose of this paper is to apply 
well known constructions in \emph{continuum theory} 
to find a space where 
the limit of $f^{-n}(A^s)$ and $f^n(A^u)$, as $n\to+\infty$,
are the stable and the unstable foliations (two different points) interpreted in this space.
If $X$ is a compact metric space then we denote by $\continua(X)$ the set of compact connected subsets of $X$.
% As noticed by Kato \cite{Ka2} for every non-trivial continuum $C\in\continua(X)$ 
% it holds that $\diam(f^n(C))\to 0$ as $n\to+\infty$ or 
% $f^n(C)\to X$ as $n\to+\infty$. 
% A similar result holds for $n\to-\infty$. 
% Let $C\subset X$ be a small ball. In this case $f^n(C)\to X$ 
% as $n\to\pm\infty$. But $C$ stretches in the direction of the unstable foliation 
% as $n\to+\infty$ and in the direction of the stable foliation as $n\to-\infty$.
% In $\continua(X)$ with the Hausdorff metric we just see that $C$ goes to $X$ but we do not see 
% how is this approximation. 
% For this purpose 
We will consider the following metric on $\continua(X)$. 

% \begin{df}
%   \label{deDistTildeH}
%   For $A,B\in\continua(X)$ the distance $\dist_{\tilde H}(A,B)$ is the minimal value of $r>0$ satisfying: 
%  1) for every subcontinuum $C\subset A$ there is a subcontinuum $D\subset B$ such that $\dist_H(C,D)\leq r$
%  and 2) for every subcontinuum $C\subset B$ there is a subcontinuum $D\subset A$ such that $\dist_H(C,D)\leq r$.
% \end{df}
\begin{df}
  \label{deDistTildeH}
 For $A_1,A_2\in\continua(X)$ the distance $\dist_{\tilde H}(A_1,A_2)$ is the minimal value of $r>0$ satisfying: 
 for $\{i,j\}=\{1,2\}$ and every subcontinuum $C\subset A_i$ there is a subcontinuum $D\subset A_j$ 
 such that $\dist_H(C,D)\leq r$.
\end{df}

In general $(\continua(X),\dist_{\tilde H})$ is not compact, as shown in Example \ref{exaHypCirculo}.
% Let us give an example. Let $X=[0,1]\times[0,1]$. 
% Consider $$C_n=\{0\}\times[0,1]\bigcup \cup_{i=0}^n[0,1]\times\{i/n\}.$$
% It is a sequence of continua converging to $X$ in the Hausdorff metric. 
% But, considering $\dist_{\tilde H}$ the sequence $C_n$ has no limit point in $\continua(X)$.
In this paper we consider a compactification $\tilde \continua(X)$ of $(\continua(X),\dist_{\tilde H})$ via \emph{hypercontinua}, 
that is, continua in $(\continua(X),\dist_H)$. 
We will see that some hypercontinua are naturally associated to foliations.
In this way we will be able to see how a continuum is stretching in the \emph{direction} of the corresponding \emph{foliation}. 
\textcolor{black}{These }ideas will be extended for cw-expansive homeomorphisms on a compact metric space. 
Recall from \cites{Ka93,Ka2} that a homeomorphism $f\colon X\to X$ is \emph{cw-expansive} 
if there is $\delta>0$ such that if $A\subset X$ is connected with more than one point 
then there is $n\in\Z$ such that $\diam(f^n(A))>\delta$.

This paper is organized as follows. 
In Section \ref{secHypercontinua} well know facts from continuum theory are recalled and examples are given. 
The reader should consult \cite{IN} for more on these topics.
In Section \ref{secPA} we apply the constructions obtained in the 
study of pseudo-Anosov maps from a viewpoint 
of continuum theory. 
In Theorems \ref{thmPApos} and \ref{thmPApos2} we show that negative iterates
of stable arcs converges to the hypercontinuum associated to the stable foliation.
In Section \ref{secHypercw} some properties of pseudo-Anosov maps are generalized 
to arbitrary cw-expansive homeomorphisms of 
a continuum.

\section{Hypercontinua}
\label{secHypercontinua}
Let $(X,\dist)$ be a continuum, i.e., a compact connected metric space.
Denote by $\continua(X)$ the set of continua contained in $X$, that is 
\[
  \continua(X)=\{A\subset X: A \text{ is compact and connected}\}.
\]
In $\continua(X)$ the Hausdorff metric is defined as
\begin{equation}
  \label{eqHdist}
    \dist_H(A_1,A_2)=\inf\{\epsilon>0:A_1\subset B_\epsilon(A_2)\text{ and } A_2\subset B_\epsilon(A_1)\},
\end{equation}
where $B_\epsilon(x)$ is the open ball of radius $\epsilon$ and centered at $x\in X$ and 
if $A\subset X$ then $B_\epsilon(A)=\cup_{x\in A}B_\epsilon(x)$.
It is known that $(\continua (X),\dist_H)$ is compact and arc-connected \cite{IN}, given that $X$ is compact and connected.
The space $\continua(X)$ is usually called the \emph{hyperspace} of $X$. 

Now we will consider a second level of abstraction.
Denote by $\continua^2(X)=\continua(\continua(X))$. 
In $\continua^2(X)$ we consider the Hausdorff distance $\dist_{H^2}$ induced 
by the metric $\dist_H$ of $\continua(X)$. 
We have again that $(\continua^2(X),\dist_{H^2})$ is compact and arc-connected.
The elements of $\continua^2(X)$ will be called \emph{hypercontinua}. 
% , continua of the hyperspace $\continua(X)$.

Each element $A\in\continua(X)$ can be seen in $\continua^2(X)$ in two natual ways. 
First, as a singleton $\{A\}\in\continua^2(X)$. 
Second, the way that will be explored in this paper, 
we can consider $\continua(A)\in\continua^2(X)$, i.e., the set formed by $A$ and its subcontinua.
Let us plainly define a map representing this idea.
Define $i\colon \continua(X)\to \continua^2(X)$ as\footnote{In fact, $i(A)=\continua(A)$.} 
\[
i(A)=\{A'\in\continua(X):A'\subset A\}.
\]
We say that $i(A)\in\continua^2(X)$ is the hypercontinuum associated to the continuum $A\in\continua(X)$.

Now we will state some properties of the map $i$.
% Define $u\colon \continua^2(X)\to\continua(X)$ by 
% $u(\tilde A)$ as the union of the sets belonging to $\tilde A$. 
% This means $u(\tilde A)=\cup_{A\in\tilde A} A$. 
% It is easy to see that $u\circ i$ is the identity of $\continua(X)$.
% Therefore, 
It is easy to prove that $i\colon \continua(X)\to\continua^2(X)$ is injective. 
It is known that the map $i$ may not be 
continuous.\footnote{In \cite{IN}, a space $X$ whose induced map $i$ is continuous is called $C^*$-\emph{smooth}. 
There, our map $i$ is called $C^*$.}
Let us give an example showing this statement.
\begin{exa}
Consider the square $X=[0,1]\times[0,1]$ with the Euclidean metric. 
Take $A$ and $A_r$ in $\continua(X)$ as in Figure \ref{figNoConti}.
\begin{figure}[h]
\center\includegraphics[scale=.7]{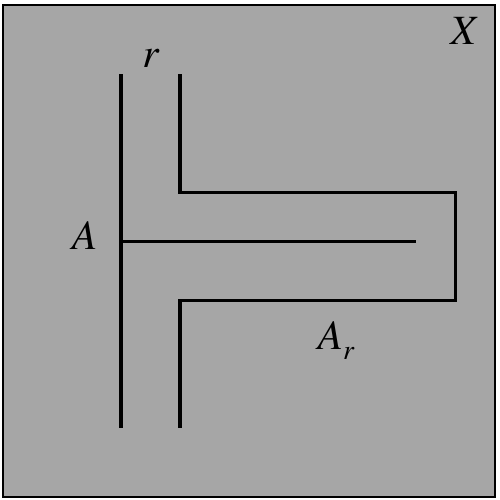}
\caption{Subcontinua of the square $X$. The parameter $r$ denotes the width of the corridor determined by $A$ and $A_r$.}
\label{figNoConti}
\end{figure}
We have that $\dist_H(A_r,A)\to 0$ as $r\to 0$.
Let us show that the distance between $i(A)$ and $i(A_r)$ is bounded away from zero.
In $\continua^2(X)$ we consider the metric $\dist_{H^2}$ induced by $\dist_H$ in $\continua(X)$.
We have that $\dist_{H^2}(i(A),i(A_r))<\epsilon$ if and only if 
for each element of $i(A)$ there is an element of $i(A_r)$ that is $\epsilon$-close,
in the metric $\dist_H$,
to the element of $i(A)$, and conversely.
This follows by the definition of the Hausdorff metric for the hypercontinua $i(A)$ and $i(A_r)$.
Now notice that the set $A$ is formed by two lines, call $A_v$ the vertical one. 
Then, $i(A)$ and $i(A_r)$ are far because there is no subcontinuum of $A_r$ being close to $A_v$, independently of how small $r$ is.
This proves that $i$ is not continuous. 
\end{exa}

% The given example suggests to consider another metric for subcontinua of $X$.
Via the injective map $i$ we define the \emph{second order Hausdorff metric} on $\continua(X)$ as 
$$\dist_{\tilde H}(A_1,A_2)=\dist_{H^2}(i(A_1),i(A_2)).$$
As the reader can check, this is exactly the metric given in Definition \ref{deDistTildeH}.
\begin{prop}
It holds that $$\dist_H(A_1,A_2)\leq\dist_{\tilde H}(A_1,A_2)$$  
for all $A_1,A_2\in\continua(X)$.
\end{prop}
\begin{proof}
  It is direct from the definitions. However we give the details. 
  Given $\epsilon>\dist_{\tilde H}(A_1,A_2)$ we will show that 
  $\dist_H(A_1,A_2)<\epsilon$. 
  If $\dist_{\tilde H}(A_1,A_2)<\epsilon$ and $x\in A_1$ (similar for $A_2$) 
  then there is a continuum $A\subset A_2$ such that 
  $\dist_H(\{x\},A)<\epsilon$. 
  Then, for every $a\in A$ we have that $\dist(x,a)<\epsilon$. 
  Consequently, $\dist_H(A_1,A_2)<\epsilon$.
\end{proof}

Define 
\[
  \tilde \continua(X)=\continua(X)\cup [\clos(i(\continua(X)))\setminus i(\continua(X))].
\]
The closure is taken in $(\continua^2(X),\dist_{H^2})$.
The metric $\dist_{\tilde H}$ is extended to $\tilde \continua(X)$.\footnote{The reader can consider 
$\tilde\continua(X)=\clos(i(\continua(X))$. However, we wish to think of $\tilde\continua(X)$ as a set extending
$\continua(X)$.} 
The space $(\tilde\continua(X),\dist_{\tilde H})$ is a compactification 
of $(\continua(X),\dist_{\tilde H})$ by adding some limit hypercontinua.
Let us give some examples to illustrate these concepts (see \cite{IN} for more on this subject). 

\begin{exa}[Hypercontinua of the interval]
Consider the interval $X=[0,1]$ with its usual metric and define the triangle 
$T=\{(a,b)\in [0,1]\times [0,1]:a\leq b\}$. 
The continua in $X$ are intervals $[a,b]\subset [0,1]$.
Therefore we can identify the hyperspace $\continua([0,1])$ with $T$ using 
the map $[a,b]\in \continua([0,1])\mapsto (a,b)\in T$.\footnote{In this case the hyperspace is two-dimensional, a triangle. 
It is known that if $X$ has topological dimension greater than 1 then the topological dimension of the hyperspace is infinite \cite{IN}.}
In this case the map $i$ is continuous. 
% The trick of Figure \ref{figNoConti} cannot be done in an arc. 
Then, $\continua([0,1])$ is compact with the metric $\dist_{\tilde H}$ 
and $\tilde\continua([0,1])=\continua([0,1])$. 
In this case no hypercontinuum is needed in the compactification.
\end{exa}

The situation is different for the circle. 

\begin{figure}[h]
  \center\includegraphics[scale=.5]{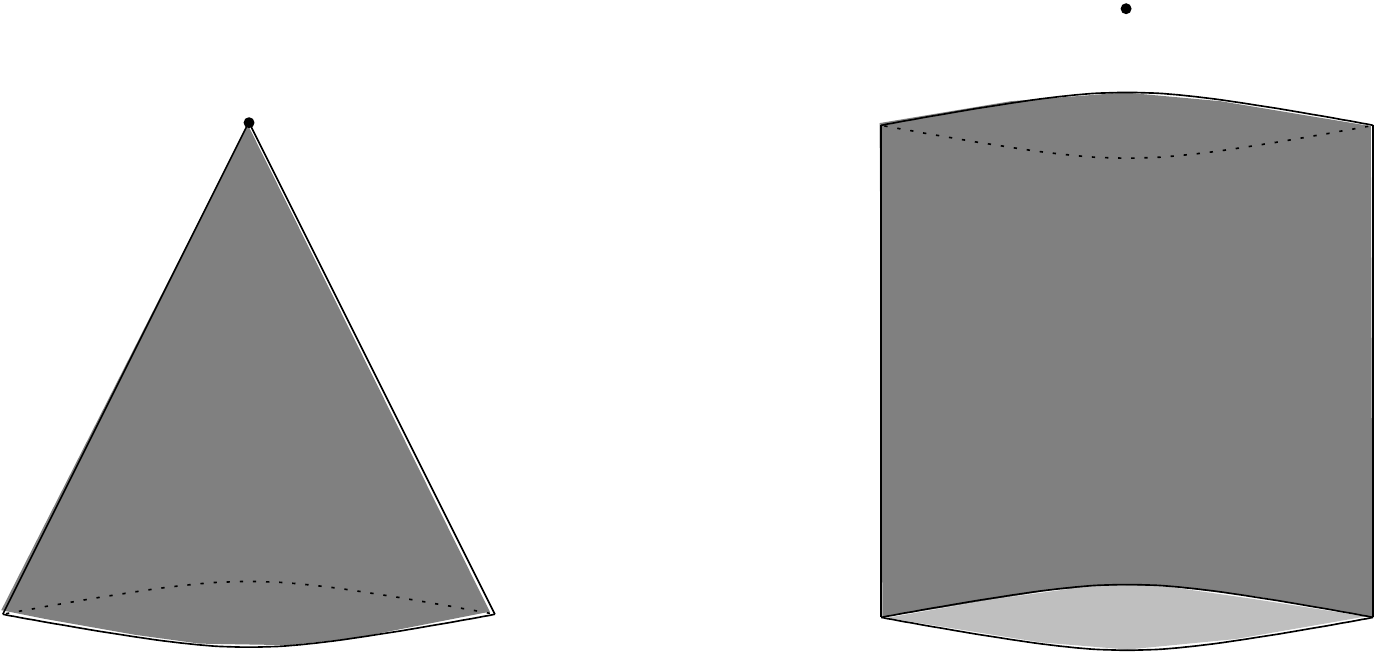}
  \caption{Left: the hyperspace $\continua(S^1)$ of the circle. 
  Right: the space $\tilde\continua(S^1)$, notice the singleton component above the cylinder.}
  \label{figHyperUnidim}
\end{figure}

\begin{exa}[Hypercontinua of the circle]
\label{exaHypCirculo}
Let $X=S^1$ be the circle $\R/\Z$ with its usual metric and orientation. 
The subcontinua of $X$ can be parameterized by its starting point and its length.
Length zero means singleton and length 1 is $X$ independently of the starting point.
Therefore, $\continua(S^1)$ can be identified with a cone (topologically a disc) 
obtained by collapsing a boundary of the cylinder $X\times [0,1]$.
Let us now obtain a model for $\tilde\continua(S^1)$. 
Let us first show that the map $i$ is not continuous. 
Consider a sequence of arcs $A_n\subset X$ with a common starting point $s\in X$ and length $1-1/n$. 
We have that the limit of $i(A_n)$ in $\tilde\continua(X)$ is the \emph{marked hypercircle} 
\[
S^1_s = \{A\in\continua(S^1): s\text{ is not in the interior of }A\}.
\]
In the Hausdorff metric, the sequence $A_n$ converges to $S^1$, but $i(A_n)$ does not converge to 
$i(S^1)$. Then, $i$ is not continuous. 
In this case, the space $(\continua(S^1),\dist_{\tilde H})$ is compactified by adding the marked hypercircles. 
In this way we obtain that $\tilde\continua(S^1)$ has two components. 
One component is a cylinder that is the union of the arcs of length less than 1 compactified with the marked hypercircles. 
The other component is the singleton $\{S^1\}$. In Figure \ref{figHyperUnidim} the spaces are illustrated.
\end{exa}

% \aca Creo que para $\dim(M)\geq 2$ siempre da conexo.

\section{Pseudo-Anosov maps}
\label{secPA}

Let $f\colon S\to S$ be a pseudo-Anosov map
of the compact orientable surface $S$ with empty boundary. 
Denote by $F^s$ and $F^u$ the stable and the unstable singular foliations of $f$, respectivelly. 
We assume that $F^s$ and $F^u$ have a finite number of singularities, being each one an 
$n$-prong with $n\geq 1$.
On the surface we consider a flat metric with conical points at the singularities of the foliations. 
It is well known that in this case the curvature is concentrated at the singularities. 
See, for example, \cite{MS} for more on this topic. 
For an $n$-prong singularity $p$ the curvature at $p$ equals $(2-n)\pi$.
The Gauss-Bonnet Theorem is true in this setting, as explained in \cite{Sch}.
\subsection{Non-negative curvature}
\label{secNNcurv}
In this section we will suppose that the curvature is non-negative at each singularity. 
The general case is considered in the next section. 
The non-negative curvature condition implies that the surface is the sphere or the torus. 
In the case of the torus there are no singular points and we will assume that $f$ is an Anosov 
diffeomorphism on a flat torus. 
In the case of the sphere there must be 4 singularities each one with curvature $\pi$, this follows by Gauss-Bonnet Theorem.
A singularity of curvature $\pi$ is usually called \emph{1-prong}.
% We consider the possibility of $f$ to have 1-prong singularities. 
% In the literature they are called generalized pseudo-Anosov. 
% It is easy to see that they are cw-expansive, see \cite{Ka2}*{Proposition 2.4}. 
\begin{rmk}
In the case of the sphere we can consider
the double branched covering of $f$ in order to obtain an Anosov torus diffeomorphism.
\end{rmk}

Now we will associate a hypercontinuum to the stable and unstable singular foliations of $f$.
We start defining
\[
\continua(F^\sigma)=\{A\in \continua(S):A\hbox{ is contained in a leaf of } F^\sigma\}
\]
for $\sigma=s,u$.
Notice that the elements of $\continua(F^\sigma)$ are arcs or singletons (trivial arcs), 
this is because we are not considering $n$-prong singularities with $n>2$.
\begin{prop}
\label{propFsConn}
The sets $\continua(F^\sigma)$, $\sigma=s,u$, are connected in 
$(\continua(S),\dist_H)$.  
\end{prop}

\begin{proof}
Take two stable arcs $A_1,A_2\in \continua(F^\sigma)$. 
Suppose that $\gamma_i\colon [0,1]\to S$ are continuous curves, $i=1,2$, 
such that $\gamma_i([0,1])=A_i$. 
Define $\alpha_j\colon[0,1]\to \continua(F^\sigma)$, for $j=1,2,3$,
such that $\alpha_1(t)=\gamma_1([0,1-t])$, $\alpha_3(t)=\gamma_2([0,t])$ 
and $\alpha_2(t)=\{\beta(t)\}$ where $\beta\colon[0,1]\to S$ is a continuous curve from $\gamma_1(0)$ 
to $\gamma_2(0)$. 
The concatenation of $\alpha_1,\alpha_2$ and $\alpha_3$ is a continuous 
curve from $A_1$ to $A_2$ inside $\continua(F^\sigma)$. 
This proves that $\continua(F^\sigma)$ is arc-connected, in particular connected.
% Notice that $F_1(S)\subset F^\sigma$ and that $F_1(S)$ is arc connected. 
% Therefore, it is sufficient to consider a curve reducing $A\in F^\sigma$ to a point of $A$ inside the leaf of $A$.
\end{proof}

Let us show that $\continua(F^\sigma)$ is not closed in $\continua(X)$.
Denote by $\len(A)$ the length of an arc in $\continua(F^\sigma)$.
Take a sequence $A_n\in \continua(F^\sigma)$ such that $\len(A_n)=n$. 
It is easy to see that $A_n$ has no convergent subsequence in $\continua(F^\sigma)$. 
Therefore, it is not closed.
Define 
$$\tilde F^\sigma=\clos(\continua(F^\sigma))$$
where the closure is taken in $(\continua(S),\dist_H)$.
From Proposition \ref{propFsConn} 
and the fact that the closure of a connected set is connected,
we know that $\tilde F^\sigma\in\continua^2(S)$, i.e., $\tilde F^\sigma$ is a hypercontinuum.

\begin{prop}
It holds that $\tilde F^\sigma=\continua(F^\sigma)\cup \{S\}$ for $\sigma=s,u$.  
\end{prop}

\begin{proof}
Let $A_n$ be a sequence in $\continua(F^\sigma)$. 
If $A_n$ is convergent $\continua(S)$, but its limit is not in $\continua(F^\sigma)$ then 
$\len(A_n)\to\infty$. 
In this case it is easy to prove that $\lim_{n\to +\infty}\dist_H(A_n,S)=0$,
which finishes the proof.
\end{proof}

Since $\continua(F^\sigma)\subset \continua(S)$ we have that 
$\tilde F^\sigma\in\tilde\continua(S)$.
In the next theorem we will show that 
$f^{-n}(A^s)$ converges to $\tilde F^s$ in the metric $\dist_{\tilde H}$ 
if $A^s$ is a non-trivial stable arc. 
For its proof we need the following simple lemma.

\begin{lem}
\label{lemArcSt}
  If $x,y\in S$ and $A_x$ is a stable arc starting at $x$ 
  then there is a stable arc $A_y$ starting at $y$ such that $\dist_H(A_x,A_y)\leq\dist(x,y)$ and 
  $\len(A_y)\leq\len(A_x)$.
\end{lem}

\begin{proof}
  In the case of the torus we can consider its Lie group structure. 
  Therefore, the stable arc $A_y$ is obtained as $A_x+y-x$, a translation of $A_x$. 
  We obtain $\len(A_y)=\len(A_x)$.
  
  In the case of the sphere we can consider its double branched covering $F$ on the torus. 
  Then, $A_y$ is obtained as the projection of the corresponding arc in the torus, as explained before. 
  Notice that if a stable arc contains a singularity, then its projection may be a shorter arc.
\end{proof}

\begin{thm}
\label{thmPApos}
If $f\colon S\to S$ is a pseudo-Anosov map with non-negative curvature then 
\[\lim_{n\to +\infty}\dist_{\tilde H}(f^{-n}(A^s),\tilde F^s)=0\]
for every non-trivial stable arc $A^s\in C(F^s)$.
Similarly for $A^u\in C(F^u)$.
\end{thm}

% We give two proofs of Theorem \ref{thmPApos}. 
% The second one is shorter but the first one is 
% straightforward. 

\begin{proof}
% [First proof of Theorem \ref{thmPApos}]
Let $A^s$ be a non-trivial stable arc. 
Recall that
$$\dist_{\tilde H}(f^n(A^s),\tilde F^s)=\dist_{H^2}(i(f^n(A^s)),\tilde F^s)$$
and notice that $i(f^n(A^s))\subset\tilde F^s$ for all $n\in\Z$.
We will show that
for all $\epsilon>0$ there is $n_0$ such that if $n>n_0$ and 
$A_1$ is a stable arc then there is 
a stable arc $A_2\subset f^{-n}(A^s)$ 
such that $\dist_H(A_1,A_2)<\epsilon$.
Given $\epsilon>0$, take $l>0$ such that 
if $A$ is a stable arc and $\len(A)\geq l$ then $\dist_H(A,S)<\epsilon/2$. 
Given a stable arc $A^s$ consider $n_0$ such that $$\len(f^{-n}(A^s))>3l$$ for all $n\geq n_0$. 
Consider an arbitrary stable arc $A_1$. 
If $\len(A_1)\geq l$ then 
$$\dist_H(A_1,f^{-n}(A^s))\leq \dist_H(A_1,S)+\dist_H(S,f^{-n}(A^s))\leq\epsilon/2+\epsilon/2=\epsilon.$$
In this case, the stable arc $A_2$ mentioned above is $A_2=f^{-n}(A^s)$. 
Now assume that $\len(A_1)<l$. 
Call $x$ an extreme point of $A_1$. 
Since $\len(f^{-n}(A^s))>3l$ we can take three disjoint sub-arcs, each one having length 
$l$. Denote by $A^s_*$ the one in the middle. 
Take $y\in A^s_*$ such that $\dist(x,y)<\epsilon$.
From Lemma \ref{lemArcSt} below we know that there is a stable arc $A_y$ such that 
$\len(A_y)\leq l$ and $\dist_H(A_y,A_1)<\epsilon$. 
Since $y\in A^s_*$ and $\len(A_y)\leq l$ we have that $A_y\subset f^{-n}(A^s)$. 
Therefore, taking $A_2=A_y$ the proof ends.
\end{proof}

% \begin{proof}[Second proof of Theorem \ref{thmPApos}]
% For $l>0$ define 
% \[
%   \continua_l(F^s)=\{A\in \continua(F^s):\len(A)\leq l\}.
% \]
% It is easy to prove that $\continua_l(F^s)$ is a closed subset of $\continua(S)$ with the metric $\dist_{\tilde H}$.
% \end{proof}

\begin{rmk}
\label{rmkPApos}
Let us explain why Theorem \ref{thmPApos} is not true for a pseudo-Anosov map 
with negative curvature. 
Suppose that $f$ has an $n$-prong singular point $p\in S$ with $n>2$. 
In this case, the local stable set $W^s_\delta(p)$ of $p$ is the union 
of $n$ arcs starting at $p$. 
We can see that no non-trivial stable arc $A^s$ can be arbitrarily close to 
$W^s_\delta(p)$. 
Since $W^s_\delta(p)\in \tilde F^s$, we have that $\dist_{\tilde H}(f^{-n}(A^s),\tilde F^s)$ is bounded away from zero.
Also, note that Lemma \ref{lemArcSt} is not true if there are $n$-prongs with $n\geq 3$. 
In the next section we explain how to adapt Theorem \ref{thmPApos} for pseudo-Anosov maps 
with singularities of negative curvature.
\end{rmk}

\subsection{Arbitrary curvature}

Let $f\colon S\to S$ be a pseudo-Anosov map
of a compact surface $S$. 
As before, denote by $F^s$ and $F^u$ the stable and the unstable singular foliations of $f$, respectivelly. 
In this section we consider $n$-prong singularities with any $n=1,3,4,\dots$.
In relation with Remark \ref{rmkPApos} we introduce the following concept.

\begin{df}We say that an arc of the foliation $F^\sigma$ is \emph{accessible} if it is contained in 
a closed product box of the foliations.
\end{df}
\begin{rmk}
  For a generalized pseudo-Anosov with non-negative curvature, as in the previous section, 
  every stable continuum is an accesible arc. 
  The local stable set of a 3-prong singularity is a non-accessible stable arc, because it is not even an arc. 
  For an $n$-prong singularity $p$ with $n>3$ there are stable arcs that are not accessible, take the union of 
  two stable arcs starting at $p$ but in non-consecutive sectors.
\end{rmk}

Define
\[
\continua_a(F^\sigma)=\{A\in \continua(S):A\hbox{ is an accessible arc of } F^\sigma\}
\]
for $\sigma=s,u$.
\begin{prop}
\label{propFsConn2}
The sets $\continua_a(F^\sigma)$, $\sigma=s,u$, are connected in 
$(\continua(S),\dist_H)$.  
\end{prop}

\begin{proof}
  Is similar to the proof of Proposition \ref{propFsConn}. 
  Note that singletons are accessible stable (and unstable) arcs.
\end{proof}

Define the hypercontinuum
$\tilde F_a^\sigma=\clos(\continua_a(F^\sigma)).$
As in the previous section, 
it holds that $\tilde F_a^\sigma=\continua_a(F^\sigma)\cup \{S\}$.  
The arguments for the following result are similar to those in the proof of Theorem \ref{thmPApos}. 
However, Lemma \ref{lemArcSt} is not true in the present case and some care is needed.

\begin{thm}
\label{thmPApos2}
If $f$ is a pseudo-Anosov map then 
\[\lim_{n\to +\infty}\dist_{\tilde H}(f^{-n}(A^s),\tilde F_a^s)=0\]
for every non-trivial accesible stable arc $A^s\in C_a(F^s)$.
\end{thm}

\begin{proof}
Fix $\epsilon>0$ given. 
Take $l_0>0$ such that if $A^s$ is a stable arc with $\len(A^s)\geq l_0$ then 
$\dist_H(A^s,S)<\epsilon/2$. 
Consider a family of closed product boxes $R_1,\dots,R_k$ of the foliations 
% and assume that $S=\cup_{j=1}^k R_j$ and 
such that for every stable arc $A^s$ with $\len(A^s)\leq l_0$ 
there is a box $R_j$ containing $A^s$. 
This boxes can be obtained from negative iterates of an arbitrary cover of product boxes.
Also assume that for every stable arc $A^s$ contained in some box $R_j$ and every point $x\in R_j$ 
there is a stable arc $A'$ such that $x\in A'\subset R_j$ and $\dist_H(A',A^s)<\epsilon$.
Take $l_1>0$ such that every stable arc of length $l_1$ cuts every box $R_j$. 
Also assume that no stable arc contained in a box $R_j$ is longer than $l_1$.

Take $A_1,A_2$ stable arcs and assume that $\len(A_2)>\max\{3l_1,l_0\}$. 
We will find a stable arc $A_3\subset A_2$ such that $\dist_H(A_1,A_3)<\epsilon$. 
If $\len(A_1)\geq l_0$ then 
\[
  \dist_H(A_1,A_2)\leq \dist_H(A_1,S)+\dist(S,A_2)< \epsilon/2+\epsilon/2=\epsilon.
\]
If $\len(A_1)<l_0$ then $A_1$ is contained in a box, say $R_{j^*}$.
Divide $A_2$ in three arcs of equal length and take the one in the middle $A_2^*$ 
with $\len(A_2^*)>l_1$. 
We know that $A_2^*$ cuts $R_{j^*}$. Take $x\in A_2^*\cap R_{j^*}$. 
The maximal stable arc of $x$ in the box $R_{j^*}$ contains a stable arc $A_3$ such that 
$\dist_H(A_3,A_1)<\epsilon$. By construction, $A_3\subset A_2$ and the proof ends.
\end{proof}

An analogous result holds for an unstable arc $A^u\in C_a(F^u)$.
It is easy to see that pseudo-Anosov maps are cw-expansive, and expansive if there are no 1-prongs. 
Also, it is known that these maps are topologically mixing and that every stable leaf is dense in the surface. 
In the next section we generalize these properties for an arbitrary cw-expansive homeomorphism on a 
continuum $X$.

\section{Hypercontinua and cw-expansivity}
\label{secHypercw}
Let $f\colon X\to X$ be a cw-expansive homeomorphism of the compact metric space $(X,\dist)$. 
% % Notice that it is usual to also denote by $f$ the map acting on subsets of $X$. 
% % That is, $f(A)=\{f(x):x\in A\}$ for $A\subset X$. 
% % This could be confusing, therefore, we introduce a new notation.
% Define 
% $\tilde f\colon \continua(X)\to \continua(X)$ as 
% $$\tilde f(A)=\{f(x):x\in A\}.$$ 
% We can entend $\tilde f$ to the homeomorphism 
% \[
% \tilde f\colon \tilde \continua(X)\to \tilde \continua(X).
% \]
% In fact, it is conjugate to the homeomorphism $\tilde{\tilde f}$ of $\clos(i(\continua(X)))$ induced by $f$.
We say that $A\in\continua(X)$ is a \emph{stable continuum} if 
$\diam(f^n(A))\to 0$ as $n\to+\infty$. 
An \emph{unstable continuum} is a stable continuum for $f^{-1}$.
Following \cite{Ka2} we define 
\[
  W^s_f=\left\{A\in\continua(X):\lim_{n\to+\infty}\diam(f^n(A))= 0\right\}
\]
and $W^u_f=W^s_{f^{-1}}$.

\begin{rmk}
  The sets $W^\sigma_f$, $\sigma=s,u$, will play the role of $\continua(F^\sigma)$ in Section \ref{secNNcurv}.
\end{rmk}

For an arbitrary homeomorphism $g\colon Y\to Y$ of a metric space and $y\in Y$ 
we define, as usual,  the $\omega$-\emph{limit set} $\omega_g(y)$ as the set of the points 
$x\in Y$ for which there is a sequence of integers $n_k\to+\infty$ such that 
$g^{n_k}(y)\to x$. 
The $\alpha$-\emph{limit set} is defined as $\alpha_g(y)=\omega_{g^{-1}}(y)$.
% Following \cite{Ka2} we define the sets of stable and unstable non-trivial continua as
% \[
%   W^s=\{A\in \continua(X)\setminus F_1(X):\diam(\tilde f^n(A))\to 0\hbox{ as } n\to+\infty\},
% \]
% \[
%   W^u=\{A\in \continua(X)\setminus F_1(X):\diam(\tilde f^n(A))\to 0\hbox{ as } n\to-\infty\}.
% \]
% Define $F_1(X)=\{\{x\}:x\in X\}$ as the space of singletons.
% Recall that $f$ is \emph{cw-expansive} if there is $\expc>0$ such that if 
% $\diam(\tilde f^n(A))\leq\expc$ for all $n\in\Z$ then $A\in F_1(X)$.
% By \cite{Ka2} we know that if $f$ is cw-expansive then $W^s\cap W^u=F_1(X)$.
% Define 
% \[
%   \omega(W^u)=\{A\in\tilde\continua(X):\exists n_k\to+\infty, B\in W^u\text{ such that }\dist_{\tilde H}(\tilde f^{n_k}(B),A)\to 0\},
% \]
% \[
%   \alpha(W^s)=\{A\in\tilde\continua(X):\exists n_k\to-\infty, B\in W^s\text{ such that }\dist_{\tilde H}(\tilde f^{n_k}(B),A)\to 0\}.
% \]
% We have that these sets are non empty, compact (with respect to $\dist_{\tilde H}$) and $\tilde f$-invariant.
% 
% \begin{rmk}
% Notice that $\alpha(W^s)\cap\omega(W^u)$ can contain some singletons. 
% For example, $Y\subset X$ is $f$-minimal (i.e. every orbit in $Y$ is dense in $Y$) 
% then $F_1(Y)\subset \alpha(W^s)\cap\omega(W^u)$.  
% \end{rmk}

Define $\tilde f\colon \tilde \continua(X)\to \tilde \continua(X)$ as the natural homeomorphism induced by $f$.
The following result represents the idea that, for a cw-expansive homeomorphism (as a Pseudo-Anosov map), 
the negative iterates of a stable continuum and 
the positive iterates of an unstable continuum wind in $X$ in different directions.

\begin{thm}
If $f$ is cw-expansive then 
$\alpha_{\tilde f}(A^s)$ and $\omega_{\tilde f}(A^u)$ are disjoint 
for every non-trivial $A^s\in W^s_f$ and $A^u\in W^u_f$.
\end{thm}

\begin{proof}
By contradiction suppose that there is $A\in\alpha_{\tilde f}(A^s)\cap\omega_{\tilde f}(A^u)$.
By definition we know that there is $n_k\to+\infty$ such that $\tilde f^{n_k}(A^u)\to A$ as 
$k\to+\infty$ in the metric $\dist_{\tilde H}$. 
Also, there is $m_k\to-\infty$ such that $\tilde f^{m_k}(A^s)\to A$ as 
$k\to+\infty$. 
We know that $A\in\tilde\continua(X)$. 

First suppose that $A\in\continua(X)$. 
Given any $\epsilon>0$ we can take a subcontinuum $A'$ of $A$ with $\diam(A')<\epsilon$. 
By definition of $\dist_{\tilde H}$, for each $k\geq 1$ there is a subcontinuum $B'_k$ of 
$\tilde f^{n_k}(A^u)$ making a convergent sequence $\dist_H(B'_k, A')\to 0$. 
If $\epsilon$ is small then we obtain that $A'$ is an unstable continuum. 
In the same way we have that $A'$ is stable. 
Since $f$ is cw-expansive we conclude that $A'$ is a singleton.
This implies that $A$ is a singleton.

Assume now that $A\in\tilde\continua(X)\setminus\continua(X)$.
Arguing as in the previous paragraph, we can prove that there is $\epsilon>0$ such that 
if $A'\in A$ and $\diam(A')<\epsilon$ then $A'$ is a singleton.
Since $A$ is a $\dist_{\tilde H}$-limit of continua we have that $A$ is a singleton. 
This contradicts that $A\in\omega_{\tilde f}(A^u)$ with $A^u$ non-trivial 
and the proof ends.
\end{proof}

% Define 
% \begin{equation}
% \begin{array}{l}
% \tilde\alpha=\alpha(W^s\setminus F_1(X)),\\
% \tilde\omega=\omega(W^u\setminus F_1(X)).
% \end{array}  
% \end{equation}
% The previous theorem implies that $\tilde\alpha$ and $\tilde\omega$ are hypercontinua. 
\begin{exa}
Consider $f\colon X\to X$ as the solenoid attractor, as explained in for example \cite{Robinson}. 
As noticed by Kato, it is positive cw-expansive (every non-trivial continuum is expanded in positive iterates).
In this case we have that stable continua are singletons.
It is known, see \cite{IN}, that the map $i$ is continuous in this case.\footnote{In terms of \cite{IN}, the solenoid is $C^*$-smooth.}
As a consequence, for every unstable non-trivial arc $A^u$ it holds that 
$$\lim_{n\to +\infty}\dist_{\tilde H}(f^n(A^u), X)=0.$$
\end{exa}

\begin{prop}
  If $X$ is a continuum then $W^s_f$ and $W^u_f$ are connected with respect to $\dist_H$.
\end{prop}

\begin{proof}
Given a stable continuum $A^s$ it is well known that there is a continuous curve $\gamma\colon[0,1]\to \continua(X)$ 
such that $\gamma(0)=A^s$, $\gamma(1)$ is a singleton while 
$\gamma(t)\subset A^s$ for all $t\in[0,1]$, see \cite{IN}. 
Now, since the space $X$ is connected the proof follows by the argument in Proposition \ref{propFsConn}.
\end{proof}

Define 
\[
\tilde F^\sigma=\clos(W^\sigma_f)
\]
for $\sigma=s,u$. 
The closure is considered in $\continua(X)$ with $\dist_H$.
We say that a homeomorphism $f$ is \emph{topologically mixing} if 
for any non-empty open sets $U$ and $V$ of $X$, there is a natural 
number $N>0$ such that $U\cap f^n(V)\neq\emptyset$ for all $n\geq N$.

The following result is similar to \cite{Ka2}*{Theorem 3.5} for fully expansive homeomorphisms.

\begin{thm}
\label{teoCwPA} If $f$ is a cw-expansive homeomorphism of the continuum $X$ and 
%   $\card(\tilde\alpha)=\card(\tilde\omega)=1$
\begin{equation}
\label{eqPA}
\left.
\begin{array}{l}
\displaystyle\lim_{n\to -\infty}\dist_{\tilde H}(f^{n}(A^s),\tilde F^s)=0,\\
\displaystyle\lim_{n\to +\infty}\dist_{\tilde H}(f^{n}(A^u),\tilde F^u)=0
\end{array}\right\}
\text{$\forall$ non-trivial $A^\sigma\in W^\sigma_f$, $\sigma=s,u$,}
\end{equation}
  then $f$ is topologically mixing.
\end{thm}

\begin{proof}
If $X$ is a singleton the result is trivial, therefore, we assume that $\dim(X)>0$ (topological dimension). 
By \cite{Ka93}*{Proposition 2.5} we know that, changing $f$ with $f^{-1}$ if needed, 
there is a non-trivial continuum $A^s\in W^s_f$.
Then, $\lim_{n\to+\infty} \dist_{\tilde H}(f^{-n}(A^s),\tilde F^s)=0$. 
By definition, for all $\epsilon>0$ there is $n_0$ such that if $n>n_0$ 
then $\dist_{\tilde H}(f^{-n}(A^s),\tilde F^s)<\epsilon$. 
In particular, for every singleton $\{x\}$ in $\tilde F^s$ there is 
$y\in f^{-n}(A^s)\cap B_\epsilon(x)$. 
Therefore, $f^{-n}(A^s)$ converges to $X$ in $\dist_H$. 
This is true for every non-trivial $A^s\in W^s_f$ and $A^u\in W^u_f$. 

% Let us show that $f$ is topologically mixing. 
Given two open subsets $U$ and $V$ of $X$, take $n$ such that $f^{-n}(A^s)$ cuts $U$. 
Now take a non-trivial subcontinuum $A_*^s\subset A^s\cap U$. 
Since $\dist_H(f^{-j}(A_*^s,X)\to 0$, there is 
$j_0$ such that for all $j>j_0$ it holds that $f^j(A_*^s)\cap V\neq\emptyset$.
This proves that $f$ is topologically mixing.
\end{proof}

\begin{rmk}
\label{rmkEstGr}
For a cw-expansive homeomorphism of a continuum $X$ satisfying (\ref{eqPA}) we have: 
  \begin{enumerate}
    \item for every non-trivial stable continuum $A^s$, $\dist_H(f^{-n}(A^s), X)\to 0$ as $n\to +\infty$ (analogous for an 
    non-trivial unstable continuum), 
    this was proved in the first paragraph of the proof of Theorem \ref{teoCwPA},
    \item if there is a non-trivial stable continuum then there is $\delta>0$ such that 
    for each $x\in X$ there is a stable continuum $A^s$ such that $\diam(A^s)=\delta$ and $x\in A^s$, 
    this follows easily from the previous property.
  \end{enumerate}

\end{rmk}

For $x\in X$ define
\[
  CW^\sigma_f(x)=\{y\in X:\exists A\in W^\sigma_f\text{ such that } x,y\in A\}
\]
for $\sigma=s,u$.

\begin{prop}
  If $X$ is a continuum and $f$ is a cw-expansive homeomorphism satisfying (\ref{eqPA}) 
  then $CW_f^\sigma(x)$ is dense in $X$ for all $x\in X$ and $\sigma=s,u$.
\end{prop}

\begin{proof}
From Remark \ref{rmkEstGr},
we can assume that there is $\delta>0$ such that 
  for every $y\in X$ there is a stable continuum $A$ such that $\diam(A)=\delta$ and $y\in A$. 
  Pick $x\in X$. For each $n\geq 0$ consider $A_n$ a stable continuum with $\diam(A_n)=\delta$ 
  and $f^n(x)\in A_n$. Define $A'_n=f^{-n}(A_n)$. 
  For each $n\geq 0$ we have that $A'_n\subset CW_f^s(x)$. 
  Then, $\cup_{n\geq 0}A'_n$ is dense in $X$. 
  Since $\cup_{n\geq 0}A'_n\subset CW_f^\sigma(x)$ the proof ends.
\end{proof}

As have seen, several dynamical properties of pseudo-Anosov maps can be deduced from 
the cw-expansivity and condition (\ref{eqPA}).
It would be interesting to characterize pseudo-Anosov maps among cw-expansive homeomorphisms 
of compact surfaces. 
For example: is every topologically mixing cw-expansive surface homeomorphism conjugate to a pseudo-Anosov map?

\noindent Departamento de Matemática y Estadística del Litoral, Salto-Uruguay\\
Universidad de la República\\
E-mail: artigue@unorte.edu.uy
\end{document}